\documentclass[11pt]{article}
\usepackage{amssymb,latexsym,amsfonts,verbatim,amscd}
\usepackage{amsmath}
\usepackage[all]{xy}

\newtheorem{Def}{Definition}[section]
\newtheorem{Thm}[Def]{Theorem}
\newtheorem{Lem}[Def]{Lemma}
\newtheorem{Prop}[Def]{Proposition}
\newtheorem{Cor}[Def]{Corollary}

\newtheorem{Que}[Def]{Question}

\newtheorem{Fac}[Def]{Fact}
\newtheorem{Ex}[Def]{Example}

\font\nat msbm10 scaled\magstephalf
\def\N{\hbox{\nat\char78}}

\def\R{\hbox{\nat\char82}}
\font\mata=msam10 
\def\restr{\mbox{\mata\char22}}
\def\telos{\hfill$\dashv$}

\def\restr{\mbox{\mata\char22}}
\def\int{\hbox{\rm int}}

\begin{document}
\sloppy

\title{Asymptotic  typicality degrees of  properties  over finite structures}
\author{Athanassios Tzouvaras}

\date{}
\maketitle

\begin{center}
Department  of Mathematics \\
Aristotle University of Thessaloniki \\
541 24 Thessaloniki, Greece. \\
e-mail: \verb"tzouvara@math.auth.gr"
\end{center}

\begin{abstract}
In previous work we defined and studied a notion of typicality,  originated with B. Russell, for properties and objects in the context of general infinite  first-order structures. In this paper we consider this  notion in the context of finite structures. In particular  we define the typicality degree of a property $\phi(x)$ over finite $L$-structures, for a language $L$, as the limit of the probability of $\phi(x)$ to be typical in an arbitrary $L$-structure ${\cal M}$  of cardinality $n$, when $n$ goes to infinity.  This poses the question whether the  0-1 law holds for typicality degrees for certain kinds of languages. One of the results of the paper is that, in contrast to the classical well-known fact that the 0-1 law holds for the  sentences of every relational language,  the 0-1 law fails  for  degrees of properties of  relational languages containing unary predicates. On the other hand it is shown that the 0-1 law holds for  degrees of some basic properties of graphs, and this gives rise to the conjecture  that the 0-1 law holds for relational languages without unary predicates. Another theme is the ``neutrality'' degree of a property $\phi(x)$ ( i.e., the fraction of $L$-structures in which neither $\phi$ nor $\neg \phi$ is typical), and in particular the ``regular'' properties (i.e., those with limit neutrality degree $0$). All properties we dealt with, either  of a relational or a  functional language, are shown to be regular,  but the question whether  {\em every} such property is regular is open.
\end{abstract}

{\em Mathematics Subject Classification (2010)}: 03C98, 03D78

\vskip 0.2in

{\em Keywords:}  Russell's notion of typicality, typical property, regular property,  0-1 law, finite graph.

\section{Typicality \`{a} la Russell}
In \cite{Tz20} we set out an investigation of a notion of typicality which is originated with B. Russell. Specifically in  \cite[p. 89]{Ru95}, Russell defines a   typical Englishman to be one ``who possesses all  the properties possessed by a majority of Englishmen.'' The notion seems captivating in its simplicity and naturalness, but in order to be formally defined one has to distinguish between properties of an object language and properties of the metalanguage, else typicality itself would be one of the properties we have to check about an  Englishman, and thus circularity arises. Once we make the aforementioned distinction using a first-order language $L$, given  any  $L$-structure  ${\cal M}=\langle M,\ldots\rangle$ we can define typical elements of $M$ in the spirit of Russell, provided we first define what a typical property is. Given a formula  $\phi(x)$ of $L$ without  parameters, let $\phi({{\cal M}})$  denote the extension of $\phi(x)$ in ${\cal M}$, i.e., $$\phi({{\cal M}})=\{a\in M:{\cal M}\models\phi(a)\}.$$

\begin{Def} \label{D:oldtypical}
{\em  Let ${\cal M}=\langle M,\ldots\rangle$ be an $L$-structure.   A property $\phi(x)$ of $L$ is said to be} typical over ${\cal M}$, {\em  if $|\phi({{\cal M}})|>|\neg\phi({{\cal M}})|=|M\backslash \phi({{\cal M}})|$.  Then an element $a\in M$ is said to be} typical {\em if it satisfies every typical property over ${\cal M}$.}
\end{Def}
In \cite{Tz20}, among other things,  we established the existence of typical elements  in many {\em infinite} structures. For example the standard structure of reals  (or second-order arithmetic) contains $|\R|$-many typical reals, while only  $<|\R|$-many nontypical ones. (A variant of the same notion of typicality,  adjusted  to fit to the context of set theory and generating a new inner model of ${\rm ZF}$, has appeared    in \cite{Tz22}.)

Instead,  in the present paper we are interested only  in typical {\em properties} (not in typical objects), and only over the {\em finite} structures of a (finite) language $L$. Specifically, we set out to study the {\em probabilities} for  $L$-properties  $\phi(x)$, in one free variable, {\em to be typical} over  arbitrary $L$-structures ${\cal M}$ of cardinality $n$, and  further to compute the limits of these probabilities, as $n$ tends to infinity. This study parallels the classical results of Finite Model Theory about the probabilities of {\em sentences} of $L$ to be {\em true} over finite structures and the fundamental  0-1 law about these  truth probabilities.

\section{Typicality degrees of first-order properties over finite structures}

\subsection{Asymptotic truth probabilities and 0-1 Laws}
Typicality degrees of properties over finite structures and their asymptotic behavior are, in some sense, generalizations of the truth degrees (or truth probabilities) of {\em sentences}. So we need to recall first some definitions and notation about the latter, see for example \cite[\S3]{Em95}, or the original paper \cite{Fa76}.  The terminology and notation here is mostly that of  \cite{Fa76}.

Let $L$ be a first-order language consisting of a finite set of (non-logical) relational symbols. For every $n$, let $\textbf{S}_n(L)$ be the set of $L$-structures ${\cal M}=\langle M,\ldots\rangle$ with $|M|=n$, or simply  $M=\{1,2,\ldots,n\}$. For every $L$-sentence $\phi$, let ${\rm Mod}_n(\phi)$ be the subset of  structures  of $\textbf{S}_n(L)$ which satisfy $\phi$. Let also
$$\mu_n(\phi)=\frac{|{\rm Mod}_n(\phi)|}{|\textbf{S}_n(L)|}, \ \mbox{and} \  \mu(\phi)=\lim_{n\rightarrow\infty}\mu_n(\phi),$$
if this limit  exists. Given a class $\Phi$ of $L$-sentences we say that $\Phi$ {\em satisfies the 0-1 law} if for every $\phi\in \Phi$, $\mu(\phi)=0$ or $1$. The following  is a fundamental result of Finite Model Theory. The following Theorem is independently due to Fagin \cite{Fa76} and Glebskii {\em et al.} \cite{GETAL69}.

\begin{Thm} \label{T:fagin}
{\rm (0-1 Law for FOL)} If $L$ is a first-order language with no function or constant symbols, then the set of  sentences of $L$ satisfies the 0-1 law.
\end{Thm}

Nevertheless  Theorem \ref{T:fagin} fails when $L$ contains function symbols. The following is the standard example  used to  prove  this failure (see \cite[\S 4]{Fa76} and \cite[Example 3.1.1]{Em95}).

\begin{Ex} \label{E:example1}
Let $L=\{F\}$, where $F$ is a unary function symbol. If $\phi$ is the $L$-sentence   $\phi:=\forall x(F(x)\neq x)$, then  $\mu(\phi)=e^{-1}$, thus the 0-1 law fails in general for the sentences of $L$.
\end{Ex}

{\em Proof.} Observe  that for any  $n\geq 1$, the number of structures  ${\cal M}=\langle M,f\rangle\in \textbf{S}_n(L)$ that satisfy $\phi:=(\forall x)(F(x)\neq x)$ is just the number of functions $f:M\rightarrow M$,  $|M|=n$, such that $f(x)\neq x$ for every $x\in M$. This number is $(n-1)^n$ (since $f(x)$ may take independently for each $x$,  $n-1$ possible values). On the other hand, $|\textbf{S}_n(L)|=n^n$. Therefore $\mu_n(\phi)=(1-1/n)^n$, and hence $\lim_{n\rightarrow\infty}\mu_n(\phi)=e^{-1}$. \telos

\subsection{Typicality degrees of properties}
Let us first elaborate a bit on  the general definition \ref{D:oldtypical}.  Recall that we denote by $\phi({\cal M})$ the extension of $\phi$ is ${\cal M}$, i.e., $\phi({\cal M})=\{a\in M:{\cal M}\models\phi(a)\}$. When we deal with typicality of {\em elements} of a structure, we naturally distinguish them  into just typical and non-typical, but when we deal with typicality of {\em properties,} especially over  finite structures,  we should  distinguish them into three kinds, according to the size of their extension.

\begin{Def} \label{D:special}
{\em Let ${\cal M}$ be an  $L$-structure  and let  $\phi(x)$ be  a property of $L$ (without parameters). We say that:

$\bullet$ $\phi(x)$ is} typical for ${\cal M}$, {\em if $|\phi({\cal M})|>|\neg\phi({\cal M})|$.

$\bullet$ $\phi(x)$ is} atypical for ${\cal M}$, {\em if $|\phi({\cal M})|<|\neg\phi({\cal M})|$.

$\bullet$ $\phi(x)$ is} neutral for ${\cal M}$, {\em if $|\phi({\cal M})|=|\neg\phi({\cal M})|$ (i.e., if neither $\phi(x)$ nor $\neg\phi(x)$ is typical). }
\end{Def}

The above distinction of properties is valid  for all structures, infinite and finite, but is particularly useful when dealing with finite structures.  If we  apply the preceding definition to a structure ${\cal M}$ of $\textbf{S}_n(L)$, then $\phi(x)$ is  typical, atypical and neutral for  ${\cal M}$, if and only if   $|\phi({\cal M})|>n/2$, $|\phi({\cal M})|<n/2$ and $|\phi({\cal M})|=n/2$, respectively, the latter case of course being possible only for even   $n$. Since for every $\phi$ and ${\cal M}\in \textbf{S}_n(L)$,
$$|\phi({\cal M})|<n/2 \Leftrightarrow |\neg\phi({\cal M})|>n/2,$$
that is,
$$\phi(x) \ \mbox{is atypical for ${\cal M}$}  \Leftrightarrow \neg\phi(x)  \ \mbox{is typical for ${\cal M}$},$$
to simplify terminology henceforth we shall not refer  to ``atypical $\phi(x)$'' but  to ``typical $\neg\phi(x)$'' instead.  Let us also fix  for every $L$ and $n$ the following subclasses of $\textbf{S}_n(L)$.
$$\textbf{S}_n(\phi:\mbox{typ})=\{{\cal M}\in \textbf{S}_n(L):\phi(x) \ \mbox{is typical for ${\cal M}$}\}=\{{\cal M}:|\phi({\cal M})|>n/2\},$$
$$\textbf{S}_n(\phi:\mbox{ntr})=\{{\cal M}\in \textbf{S}_n(L):\phi(x) \ \mbox{is neutral for ${\cal M}$}\}=\{{\cal M}:|\phi({\cal M})|=n/2\}.$$
The second of the above  classes exists only for even $n$, so for each property  $\phi(x)$, $\textbf{S}_n(L)$ splits as follows  for odd and even $n$:
\begin{equation} \label{E:odd}
\textbf{S}_{2n+1}(L)=\textbf{S}_{2n+1}(\phi:\mbox{typ})\cup \textbf{S}_{2n+1}(\neg\phi:\mbox{typ}),
\end{equation}
while
\begin{equation} \label{E:even}
\textbf{S}_{2n}(L)=\textbf{S}_{2n}(\phi:\mbox{typ})\cup \textbf{S}_{2n}(\neg\phi:\mbox{typ})\cup\textbf{S}_{2n}(\phi:\mbox{ntr}).
\end{equation}
Then, by analogy with the  probabilities $\mu_n(\phi)$, and asymptotic probability $\mu(\phi)=\lim_{n\rightarrow\infty}\mu_n(\phi)$ for the truth of $L$-sentences referred to in Section 2.1 above, we can naturally define the corresponding probabilities for a property $\phi(x)$ to be typical or neutral over an arbitrary structure ${\cal M}\in \textbf{S}_n(L)$. Specifically, for each  $n$, we  set
$$d_n(\phi:\mbox{typ})=\frac{|\textbf{S}_n(\phi:\mbox{typ})|}{|\textbf{S}_n(L)|}, \quad \ d_{2n}(\phi:\mbox{ntr})=\frac{|\textbf{S}_{2n}(\phi:\mbox{ntr})|}
{|\textbf{S}_{2n}(L)|},$$
and, further,
$$d(\phi:\mbox{typ})=\lim_{n\rightarrow\infty}d_n(\phi:\mbox{typ}), \quad d(\phi:\mbox{ntr})=\lim_{n\rightarrow\infty}d_{2n}(\phi:\mbox{ntr}),$$
whenever these limits exist. $d_n(\phi:\mbox{typ})$ and $d_{2n}(\phi:\mbox{ntr})$ are the {\em $n$-typicality degree} and {\em $n$-neutrality degree} of $\phi(x)$, respectively, while  $d(\phi:\mbox{typ})$ and $d(\phi:\mbox{ntr})$ are the corresponding {\em asymptotic degrees}. Here are some   first  consequences of the  definitions.

\begin{Fac} \label{F:Boolean}
(i) If $\vdash \phi(x)\rightarrow \psi(x)$, then $d_n(\phi:\mbox{typ})\leq d_n(\psi:\mbox{typ})$, for all $n\geq 1$. Therefore, if $d(\psi:\mbox{typ})=0$, then  $d(\phi:\mbox{typ})=0$ too.

(ii) If $\vdash \phi(x)\leftrightarrow \psi(x)$, then $d_n(\phi:\mbox{typ})=d_n(\psi:\mbox{typ})$, and also $d_{2n}(\phi:\mbox{ntr})=d_{2n}(\psi:\mbox{ntr})$ for all $n\geq 1$.

(iii) For all $\phi(x)$ and $n$, $d_{2n}(\phi:\mbox{ntr})=d_{2n}(\neg\phi:\mbox{ntr})$.

(iv) If $d(\phi:\mbox{typ})=a>0$ (resp. $d(\phi:\mbox{ntr})=a>0$), then the set $\{n:\textbf{S}_n(\phi:\mbox{typ})\neq \emptyset\}$ (resp. $\{n:\textbf{S}_{2n}(\phi:\mbox{ntr})\neq \emptyset\}$) is cofinite.
\end{Fac}

{\em Proof.} For (i) and (ii) just note that if $\vdash \phi(x)\rightarrow \psi(x)$, then for every structure ${\cal M}$, $\phi({\cal M})\subseteq \psi({\cal M})$, hence $|\phi({\cal M})|\leq |\psi({\cal M})|$, so
$${\cal M}\in \textbf{S}_n(\phi:\mbox{typ})\Leftrightarrow |\phi({\cal M})|>n/2\Rightarrow |\psi({\cal M})|>n/2\Leftrightarrow{\cal M}\in \textbf{S}_n(\psi:\mbox{typ}),$$
therefore $\textbf{S}_n(\phi:\mbox{typ})\subseteq \textbf{S}_n(\psi:\mbox{typ})$. Moreover $\vdash \phi(x)\leftrightarrow \psi(x)$ implies $\textbf{S}_n(\phi:\mbox{typ})=\textbf{S}_n(\psi:\mbox{typ})$ and also $\textbf{S}_{2n}(\phi:\mbox{ntr})=\textbf{S}_{2n}(\psi:\mbox{ntr})$, for every $n$.

(iii) If $|M|=2n$, then for every $\phi(x)$, obviously $|\phi({\cal M})|=n\Leftrightarrow |\neg\phi({\cal M})|=n$, so  $\textbf{S}_{2n}(\phi:\mbox{ntr})=\textbf{S}_{2n}(\neg\phi:\mbox{ntr})$.

(iv) Let $d(\phi:\mbox{typ})=a>0$. If we pick some $0<\varepsilon<a$, then clearly there is $n_0$ such that  for all $n\geq n_0$, $d_n(\phi:\mbox{typ})>a-\varepsilon$. Since $d_n(\phi:\mbox{typ})=\frac{|\textbf{S}_n(\phi:\mbox{\small typ})|}{|\textbf{S}_n(L)|}$, it follows that  for all $n\geq n_0$  $|\textbf{S}_n(\phi:\mbox{typ})|\neq 0$, and hence  $\textbf{S}_n(\phi:\mbox{typ})\neq \emptyset$. The claim for $\{n:\textbf{S}_{2n}(\phi:\mbox{ntr})\neq \emptyset\}$ is shown similarly.
\telos

\vskip 0.2in

By the definition of $\mu(\phi)$ in Section 2.1, it follows immediately that for any language $L$ and any $L$-sentence  $\phi$, if $\mu(\phi)$ exists, then so does $\mu(\neg\phi)$ and  $\mu(\neg\phi)=1-\mu(\phi)$. What about typicality degrees? Is it true that $d(\neg\phi:\mbox{typ})=1-d(\phi:\mbox{typ})$ whenever $d(\phi:\mbox{typ})$ exists, for any property $\phi(x)$? The question is eventually open and the  reason is the limit $\lim_{n\rightarrow\infty}d_{2n}(\phi:\mbox{ntr})$.  Namely, while by  (\ref{E:odd})
$$d_{2n+1}(\phi:\mbox{typ})+d_{2n+1}(\neg\phi:\mbox{typ})=1,$$
(\ref{E:even}) implies that
$$d_{2n}(\phi:\mbox{typ})+d_{2n}(\neg\phi:\mbox{typ})+d_{2n}(\phi:\mbox{ntr})=1.$$
Thus in the second case we have
$$\lim_{n\rightarrow\infty}d_{2n}(\neg\phi:\mbox{typ})=
1-\lim_{n\rightarrow\infty}d_{2n}(\phi:\mbox{typ})-\lim_{n\rightarrow\infty}d_{2n}(\phi:\mbox{ntr}),$$
and  in order to infer that   $\lim_{n\rightarrow\infty}d_{2n}(\neg\phi:\mbox{typ})=1-\lim_{n\rightarrow\infty}d_{2n}(\phi:\mbox{typ})$, we must establish  that  $\lim_{n\rightarrow\infty}d_{2n}(\phi:\mbox{ntr})=0$. We don't know if this is the case for every property $\phi(x)$ of every language. So we shall give a name to properties satisfying  this interesting and convenient condition.

\section{Regularity of properties}

\begin{Def} \label{D:regular}
{\em  A property $\phi(x)$ of $L$ is said to be} regular {\em if $d(\phi:\mbox{ntr})=0$.}
\end{Def}

\begin{Fac} \label{F:sum}
(i) $\phi(x)$ is regular if and only if  $\neg\phi(x)$ is regular.

(ii) If $d(\phi:\mbox{typ})$ exists, then so does $\lim_{n\rightarrow\infty}d_{2n+1}(\neg\phi:\mbox{typ})$ and $$\lim_{n\rightarrow\infty}d_{2n+1}(\neg\phi:\mbox{typ})=1-d(\phi:\mbox{typ}).$$

(iii) If  $\phi(x)$ is regular, then also
$$\lim_{n\rightarrow\infty}d_{2n}(\neg\phi:\mbox{typ})=1-d(\phi:\mbox{typ}).$$
and therefore
$$d(\neg\phi:\mbox{typ})=1-d(\phi:\mbox{typ}).$$
\end{Fac}

{\em Proof.}  (i) By Fact \ref{F:Boolean} (iii), for every $n$, $d_{2n}(\phi:\mbox{ntr})=d_{2n}(\neg\phi:\mbox{ntr})$, therefore $d(\phi:\mbox{ntr})=0$ if and only if  $d(\neg\phi:\mbox{ntr})=0$.

(ii) If $d(\phi:\mbox{typ})=a$, then also $\lim_{n\rightarrow\infty}d_{2n+1}(\phi:\mbox{typ})=a$, thus by (\ref{E:odd}), $\lim_{n\rightarrow\infty}d_{2n+1}(\neg\phi:\mbox{typ})=1-a$.

(iii) If $d(\phi:\mbox{typ})=a$ and $\phi(x)$ is regular, then $\lim_{n\rightarrow\infty}d_{2n}(\phi:\mbox{ntr})=0$, so  by (\ref{E:even})
$$\lim_{n\rightarrow\infty}d_{2n}(\neg\phi:\mbox{typ})=1-a-\lim_{n\rightarrow\infty}d_{2n}(\phi:\mbox{ntr})=1-a.$$
 \telos

All specific properties we treat  below are  regular. So it is natural to ask whether  every property is regular. The question is open for general languages. In the next subsection we show that it is true  for a large class of properties of the language  $L=\{U_1,\ldots,U_k\}$ which consists  of an arbitrary number of unary predicates.

\subsection{Regularity of properties of $L=\{U_1,\ldots,U_k\}$}

Let $L=\{U_1,\ldots,U_k\}$ be a language with $k$ unary predicates. For each $i\in\{1,\ldots,k\}$, let  $U_i^{1}(x)=U_i(x)$ and $U_i^{0}(x)=\neg U_i(x)$. Then given a function $e\in \{0,1\}^k$, we set   $\phi_e(x)=U^{e(1)}_1(x)\wedge U^{e(2)}_2(x)\wedge\cdots\wedge U^{e(k)}_k(x)$. The properties $\phi_e(x)$, for $e\in \{0,1\}^k$, form the $2^k$  atoms of the (syntactic) Boolean algebra ${\cal B}_{prop}$ generated by the properties $U_i(x)$, $i\in \{1,\ldots,k\}$, and any two distinct atoms $\phi_{e_1}(x)$, $\phi_{e_2}(x)$ are mutually inconsistent, i.e., $\phi_{e_1}(x)\wedge \phi_{e_2}(x)\vdash \bot$. Besides  each of  the $2^{2^k}$ elements of ${\cal B}_{prop}$ has the form   $\phi_E(x)=\bigvee_{e\in E}\phi_e(x)$, for some $E\subseteq \{0,1\}^k$.

We shall generalize the class of formulas $\phi_e(x)$ defined above, by relaxing the condition that for every $i\leq k$,  either $U_i(x)$ or $\neg U_i(x)$ must be a conjunct of $\phi_e$. Namely for any $p$ such that  $1\leq p\leq k$,  a $p$-subsequence of $\langle 1,\ldots,k\rangle$ is a  $p$-tuple $\langle i_1,\ldots,i_p\rangle$ such that  $1\leq i_1<i_2<\cdots<i_p\leq k$. Given a $p$-subsequence $\bar{s}=\langle i_1,\ldots,i_p\rangle $ and an $e\in\{0,1\}^k$, let $$\phi_{\bar{s},e}(x)=U^{e(i_1)}_{i_1}(x)\wedge\cdots\wedge  U^{e(i_p)}_{i_p}(x).$$
We refer to formulas $\phi_{\bar{s},e}(x)$ as {\em basic formulas} of $L$. The main result of this subsection is that every basic formula of $L$ is regular.

Given an $L$-structure ${\cal M}=\langle M,W_1,\ldots,W_k\rangle$, let $W_i^1=W_i$ and $W_i^0=M\backslash W_i$. For every $e\in\{0,1\}^k$, let $X_e=W^{e(1)}_1\cap W^{e(2)}_2\cap\cdots\cap W^{e(k)}_k$. Clearly the sets $X_e$, for  $e\in \{0,1\}^k$, are pairwise disjoint, but their difference from $\phi_e(x)$ is that not all of them need to be nonempty. Let  ${\cal B}_{set}\subseteq {\cal P}(M)$
be the Boolean algebra generated by the sets $W_i$, $i\in\{1,\ldots,k\}$ with at  most $2^k$ atoms.  As before  for each $E\subseteq \{0,1\}^k$, we set   $X_E=\bigcup_{e\in E}X_e$. Further for any $p$-subsequence $\bar{s}=\langle i_1,\ldots,i_p\rangle$ of $\langle 1,\ldots,k\rangle$, with $1\leq p\leq k$,  and any $e\in\{0,1\}^k$, we write $X_{\bar{s},e}=W^{e(i_1)}_{i_1}\cap\cdots\cap  W^{e(i_p)}_{i_p}$. Finally  each set $X_E$ is defined in ${\cal M}$ by the property $\phi_E(x)$, i.e., $X_E=\phi_E({\cal M})$, and  each $X_{\bar{s},e}$ is defined by $\phi_{\bar{s},e}(x)$, i.e., $\phi_{\bar{s},e}({\cal M})=X_{\bar{s},e}$.

\begin{Lem} \label{L:beginwith}
Given an $L$-structure ${\cal M}=\langle M,W_1,\ldots,W_k\rangle$ as above, the definable (without parameters) subsets of $M$ are  exactly the sets $X_E$. Therefore every $L$-property $\phi(x)$ is equivalent over ${\cal M}$ to a $\phi_e(x)$, for some $E\subseteq\{0,1\}^k$.
\end{Lem}

{\em Proof.} Towards reaching a contradiction assume that there is an $L$-property $\phi(x)$ such that  $\phi({\cal M})=A$ and  $A\neq X_E$ for every $E\subseteq \{0,1\}^k$. Since the sets $X_{e}$ form a partition of $M$, in order for a set $Y\subseteq M$ to have the property
$$(\forall e)(\forall a,b\in X_{e})(a\in Y\Leftrightarrow b\in Y),$$
it is necessary and sufficient that $Y\in {\cal B}$, where ${\cal B}$ is the Boolean algebra mentioned above, i.e., $Y=X_E$ for some $E\subseteq \{0,1\}^k$. So since by assumption $A\notin {\cal B}$ it follows that  there exists $e\in\{0,1\}^k$ such that
\begin{equation} \label{E:differ}
(\exists a,b\in X_{e})(a\in A\Leftrightarrow b\notin A),
\end{equation}
Now observe that every bijection $f:M\rightarrow M$ which preserves the sets $W_i$, i.e., such that $f[W_i]=W_i$ for every $i=1,\ldots,k$, is an automorphism of ${\cal M}=\langle M,W_1,\ldots,W_k\rangle$. By (\ref{E:differ}) we can  pick $X_{e}$ and $a,b\in X_{e}$ such that $a\in A\Leftrightarrow b\notin A$ and  take the bijection  $f_1:X_{e}\rightarrow X_{e}$ which interchanges $a$ and $b$. Let also  $id_{M\backslash X_{e}}$ be the identity on the complement of $X_{e}$ and $f=f_1\cup id_{M\backslash X_{e}}$. Then $f$ is an automorphism of ${\cal M}$ such that $f(a)=b$. Since we assumed that there is $\phi(x)$ such that $A=\phi({\cal M})$, $f$ must preserve $A$. But
$$a\in A \Leftrightarrow a\in \phi({\cal M}) \Leftrightarrow{\cal M}\models \phi(a)
\Leftrightarrow {\cal M}\models \phi(f(a)) \Leftrightarrow f(a)\in A \Leftrightarrow b\in A,$$
a contradiction.  This completes the proof. \telos

\vskip 0.2in

Let us notice here a result  which is involved in all  proofs of regularity of  properties. Whenever we try to check the regularity of a property over a structure $M$ with $2n$ elements,  we shall necessarily deal with  the number $\binom{2n}{n}$ which counts the  subsets $M$  having  half of its elements. The numbers $\binom{2n}{n}$ are known as  ``central binomial coefficients'' and several useful combinatorial facts are known about them (see e.g.  \cite{Wiki}). In particular the following upper bound is particularly helpful  and will be used below.

\begin{Fac} \label{F:central}
{\rm (\cite{Wiki})} For every $n\geq 1$,  $\binom{2n}{n}\leq \frac{4^n}{\sqrt{\pi n}}$.
\end{Fac}

For  $k\leq n$,  we shall  denote by $(n)_k$ the number of $k$-tuples of distinct elements chosen from a set of $n$ elements. It is well known that
$$(n)_k=\frac{n!}{(n-k)!}=(n-k+1)(n-k+2)\cdots n.$$
Then $(n)_1=n$, $(n)_n=n!$, $(n)_k=0$ for $k>n$. In particular we set $(n)_0=1$. The numbers $(n)_m$  are called {\em falling factorials}. Notations $n^{\underline{m}}$ and  $P(n,m)$ are often used in the bibliography instead of $(n)_m$.

Below  we shall  employ the relation  $f(n) \sim  g(n)$  of {\em asymptotic equality} between  functions $f,g:\N\rightarrow \R$, as well as that of {\em asymptotic inequality} $f(n) \lesssim  g(n)$.   $f(n) \sim  g(n)$  means, by definition,  $\lim_{n\rightarrow\infty}\frac{f(n)}{g(n)}=1$, or equivalently, $f(n)=g(n)+o(g(n))$, where $o(g(n))$ is a function such that $\lim_{n\rightarrow\infty}\frac{o(g(n))}{g(n)}=0$. $f(n) \lesssim  g(n)$ means   $f(n)\leq g(n)+o(g(n))$. The properties of $\sim$ and $\lesssim$ we shall need are the following, and are either  well-known or easily verified.

\begin{Fac} \label{F:asymlimit}
(i) The relation $\sim$ is preserved by  the usual operations, i.e.  if $f_1\sim g_1$ and $f_2\sim g_2$, then $f_1+f_2\sim g_1+g_2$, $f_1\cdot f_2\sim g_1\cdot g_2$, and $\frac{f_1}{f_2}\sim \frac{g_1}{g_2}$.

(ii) If $f(n)\sim g(n)$ and $\lim_{n\rightarrow\infty}f(n)=a\in\R$, then $\lim_{n\rightarrow\infty}g(n)=a$.

(iii) If $f(n)\sim g(n)$ and $g(n)\leq h(n)$, then $f(n)\lesssim h(n)$.

(iv) If $0\leq f(n)\lesssim g(n)$ and $\lim_{n\rightarrow\infty}g(n)=0$, then $\lim_{n\rightarrow\infty}f(n)=0$.
\end{Fac}

[For completeness we  sketch the proof of (ii). Let $f(n)\sim g(n)$ and $\lim_{n\rightarrow\infty}f(n)=a\in\R$. Then clearly both $f,g$ are bounded, and let $b$ be a bound for $g$, i.e.,  $\forall n \ |g(n)|\leq b$. We have   $f(n)=g(n)+o(g(n))$, where $\lim_{n\rightarrow\infty}\frac{o(g(n))}{g(n)}=0$. So for all $n$:

(1) $|g(n)-a|\leq |g(n)-f(n)|+|f(n)-a|=|o(g(n))|+|f(n)-a|$.

\noindent Fix some  $\varepsilon>0$. There is $n_1$ such that $\forall n\geq n_1 \ |\frac{o(g(n))}{g(n)}|\leq \frac{\varepsilon}{2b}$, hence:

(2) $\forall n\geq n_1 \ |o(g(n))|\leq \frac{\varepsilon}{2b}|g(n)|\leq \frac{\varepsilon}{2}$.

\noindent Also there is $n_2$ such that:

(3) $\forall n\geq n_2 \ |f(n)-a|<\frac{\varepsilon}{2}$.

\noindent If $n_0=\max(n_1,n_2)$, (1), (2) and (3)  yield $\forall n\geq n_0 \ |g(n)-a|<\varepsilon$.]

\vskip 0.2in

For example, for any  fixed $k$ such that $1\leq k< n$, $(n)_k$ is a polynomial in $n$ of degree $k$ with leading coefficient $1$, so
\begin{equation} \label{E:asequal}
(n)_k\ \sim \ n^k.
\end{equation}
We shall apply  relation (\ref{E:asequal}) several times  without explicit reference to that. We shall also make use of the following result.

\begin{Fac} \label{F:semi-sum}
$\sum_{k<n/2}\binom{n}{k}=\sum_{k>n/2}\binom{n}{k} \ \sim 2^{n-1}$.
\end{Fac}

{\em Proof.} The first equality is well-known. Moreover for every $n$, $\sum_{k=0}^n\binom{n}{k}=2^n$. If $n$ is odd, $\sum_{k=0}^n\binom{n}{k}=\sum_{k<n/2}\binom{n}{k}+\sum_{k>n/2}\binom{n}{k}$, so $\sum_{k<n/2}\binom{n}{k}=\frac{2^n}{2}=2^{n-1}$. If $n$ is even and $n=2m$, then $2(\sum_{k<n/2}\binom{n}{k})+\binom{2m}{m}=2^n$, that is, $\sum_{k<n/2}\binom{n}{k}=\frac{1}{2}(2^n-\binom{2m}{m})=
2^{n-1}-\frac{1}{2}\binom{2m}{m}$. Therefore
$$\frac{\sum_{k<n/2}\binom{n}{k}}{2^{n-1}}=1-\frac{1}{2^{2m}}\binom{2m}{m}.$$
So it suffices to see that $\lim_{m\rightarrow \infty}\frac{1}{2^{2m}}\binom{2m}{m}=0$. But this follows from Fact \ref{F:central}, since $\frac{1}{2^{2m}}\binom{2m}{m}\leq \frac{1}{2^{2m}}\frac{4^m}{\sqrt{\pi m}}=\frac{1}{\sqrt{\pi m}}\longrightarrow_m0$. \telos

\vskip 0.1in

\begin{Thm} \label{T:regularity}
Every  basic property $\phi_{\bar{s},e}(x)$ of $L=\{U_1,\ldots,U_k\}$ is regular. In particular every  property $\phi_e(x)$, as well as every $U_i(x)$ and $\neg U_i(x)$, for $i=1,\ldots,k$, is regular.
\end{Thm}

{\em Proof.} Let us fix the ground set $M$ of all ${\cal M}=\langle M,W_1,\ldots,W_k\rangle\in \textbf{S}_{2n}(L)$, i.e. $|M|=2n$, and fix also a basic formula $\phi_{\bar{s},e}(x)$,  for some $p$-subsequence $\bar{s}=\langle i_1,\ldots,i_p\rangle$ of  $\langle 1,\ldots,k\rangle$ and some $e\in \{0,1\}^k$. We have to compute the limit of the probability
$$d_{2n}(\phi_{\bar{s},e}:\mbox{ntr})=
\frac{|\textbf{S}_{2n}(\phi_{\bar{s},e}:\mbox{ntr})|}{|\textbf{S}_{2n}(L)|}.$$

Note that each $L$-structure ${\cal M}$ is determined by a $k$-tuple  $\langle W_1,\ldots,W_k\rangle$ of elements of ${\cal P}(M)$, rather than a $k$-element subset $\{W_1,\ldots,W_k\}$. This is because an interpretation of $L$ in ${\cal M}$ is a mapping $I:\{U_1,\ldots,U_k\}\rightarrow {\cal P}(M)$, or $I:\{1,\ldots,k\}\rightarrow {\cal P}(M)$, such that $I(i)=W_i=U_i^{\cal M}$. Each such $I$ determines a $k$-tuple  $\langle W_1,\ldots,W_k\rangle$.  To be precise,  each $W_i$ must be  different from $\emptyset$ and $M$,  but this does not affect the asymptotic behavior of the neutrality degree. Namely, by (\ref{E:asequal}),
\begin{equation} \label{E:denom}
|\textbf{S}_{2n}(L)|=(2^{2n-2})_k \ \sim \ 2^{2kn}.
\end{equation}

In order to compute $|\textbf{S}_{2n}(\phi_{\bar{s},e}:\mbox{ntr})|$ we fix temporarily a set $A\subseteq M$ such that $|A|=n$. Since  $\phi_{\bar{s},e}({\cal M})=W^{e(i_1)}_{i_1}\cap\cdots \cap W^{e(i_p)}_{i_p}$, we  set
$$Z(A)=\{\langle W_1,\ldots,W_k\rangle: W^{e(i_1)}_{i_1}\cap\cdots \cap W^{e(i_p)}_{i_p}=A\}.$$
Then clearly
\begin{equation} \label{E:cases}
|\textbf{S}_{2n}(\phi_{\bar{s},e}:\mbox{ntr})|=|Z(A)|\cdot\binom{2n}{n}.
\end{equation}
Now $|Z(A)|=|Z_1(A)|\cdot|Z_2(A)|$, where
$$Z_1(A)=\{\langle W^{e(i_1)}_{i_1},\cdots, W^{e(i_p)}_{i_p}\rangle :W^{e(i_1)}_{i_1}\cap\cdots \cap W^{e(i_p)}_{i_p}=A\},$$
and
$$Z_2(A)=\{\langle W_{j_1}^{e(j_1)},\ldots, W_{j_{k-p}}^{e(j_{k-p})}\rangle:
\{W_{j_1}^{e(j_1)},\ldots, W_{j_{k-p}}^{e(j_{k-p})}\}\subseteq$$
$$\subseteq{\cal P}(M)\backslash \{W^{e(i_1)}_{i_1},\ldots, W^{e(i_p)}_{i_p}\}\}.$$
In order to compute (or find upper bounds for)  $|Z_1(A)|$ and $|Z_2(A)|$, we must  distinguish the cases $p=1$ and  $2\leq p\leq k$.

\vskip 0.1in

{\em Case 1.} $p=1$. Without loss of generality we may assume that $\bar{s}$ is the 1-subsequence $(1)$, i.e., $\phi_{\bar{s},e}(x)$ is either $U_1(x)$ or $\neg U_1(x)$. By Fact \ref{F:sum} (i), it suffices to consider only $U_1(x)$. Then  $U_1({\cal M})=W_1$, so $Z_1(A)=\{W_1:W_1=A\}$, and hence $|Z_1(A)|=1$. Also
$Z_2(A)=\{\langle W_2,\ldots,W_k\rangle: \{W_2,\ldots,W_k\}\subseteq {\cal P}(M)\backslash \{W_1\}\}$, so, in view of  (\ref{E:asequal}),
$$|Z_2(A)|=(2^{2n}-1)_{k-1} \ \sim \ 2^{2(k-1)n}.$$
Thus also $|Z(A)|=|Z_2(A)|\ \sim \ 2^{2(k-1)n}$, and therefore, letting $A$ range over all sets of cardinality $n$,
$$\textbf{S}_{2n}(U_1(x):\mbox{ntr})\sim \binom{2n}{n}\cdot 2^{2(k-1)n}.$$
From  the last equality, (\ref{E:denom}) and Fact  \ref{F:central} we get
$$d_{2n}(U_1(x):\mbox{ntr}) \sim  \binom{2n}{n}\cdot\frac{2^{2(k-1)n}}{2^{kn}}  \sim   \binom{2n}{n}\cdot\frac{1}{2^{2n}}\leq\frac{4^n}{\sqrt{\pi n}}\cdot\frac{1}{2^{2n}}=\frac{1}{\sqrt{\pi n}}\longrightarrow_n=0.$$
Therefore  $U_1(x)$ is regular.

\vskip 0.2in

{\em Case 2.}  $2\leq p\leq k$. Fix a $p$-subsequence of $\langle 1,\ldots,k\rangle$ and an $e\in \{0,1\}^p$. Then $\phi_{\bar{s},e}({\cal M})=W_{i_l}^{e(i_1)}\cap\cdots\cap W_{i_p}^{e(i_p)}$. Fixing temporarily a set $A\subseteq M$, as before $Z_2(A)$ consists of the $(k-p)$-tuples  of ${\cal P}(M)\backslash \{W^{e(i_1)}_{i_1},\ldots, W^{e(i_p)}_{i_p}\}$, so
$$|Z_2(A)|=(2^{2n}-p)_{k-p}\sim  2^{2(k-p)n}.$$
The main difference of this case from the previous one lies in the computation of $Z_1(A)$. Observe that  $W_{i_l}^{e(i_1)}\cap\cdots\cap W_{i_p}^{e(i_p)}=A$ implies  $A\subseteq W_{i_j}^{e(i_j)}$,  for  each $j=1,\ldots,p$. For each $i_j$ we consider the cases $e(i_j)=1$ and $e(i_j)=0$.

(a) If $e(i_j)=1$, then $A\subseteq W_{i_j}$, so $W_{i_j}=A\cup Y_j$, for some  $Y_j\subseteq M\backslash A$.

(b) If $e(i_j)=0$, then  $A\subseteq M\backslash W_{i_j}$, so $M\backslash W_{i_j}=A\cup Y_j$ for some  $Y_j\subseteq M\backslash A$.

Thus in both cases for each $j=1,\ldots,p$ there are at most as many possible choices for   $W_{i_j}$ as are the choices for $Y_j\subseteq M\backslash A$, i.e.,  $2^{|M\backslash A|}=2^n$. In fact  the choices of such $Y_{i_j}$'s are not independent. They must satisfy the condition $\bigcap_{j=1}^pY_{i_j}=\emptyset$ (since otherwise $\bigcap_{j=1}^p W_{i_j}^{e(i_j)}\neq A$). Nevertheless, the number of possible $p$-tuples of $Z_1(A)$ are at most as many as the  $p$-tuples   of ${\cal P}(M\backslash A)$, i.e.,  $(2^n)_{p}\sim 2^{pn}$.  Therefore
$$|Z_1(A)|\leq (2^n)_{p}\sim 2^{pn}.$$
So
$$|Z(A)|=|Z_1(A)|\cdot|Z_2(A)|\lesssim 2^{2(k-p)n}\cdot 2^{pn}\sim 2^{(2k-p)n}.$$
Letting $A$ range over all sets of cardinality $n$, we have
$$|\textbf{S}_{2n}(\phi_{\bar{s},e}:\mbox{ntr})|\lesssim 2^{(2k-p)n}\cdot \binom{2n}{n}.$$
So, by (\ref{E:denom}) and Fact  \ref{F:central},
$$d_{2n}(\phi_{\bar{s},e}:\mbox{ntr})\lesssim  \frac{2^{(2k-p)n}\cdot\binom{2n}{n}}{2^{2kn}}\ \sim \ \frac{\binom{2n}{n}}{2^{pn}}\leq \frac{4^n}{ 2^{pn}\cdot\sqrt{\pi n}}.$$
Since $p\geq 2$, $\frac{4^n}{2^{pn}\cdot \sqrt{\pi n}}\leq\frac{4^n}{4^n\cdot \sqrt{\pi n}}=\frac{1}{\sqrt{\pi n}}$. So $d_{2n}(\phi_{\bar{s},e}:\mbox{ntr})\longrightarrow_n0$, according to Fact \ref{F:asymlimit}. This completes the proof. \telos

\vskip 0.2in

It is still open however whether {\em every} property of $L$, i.e., every $\phi_E(x)=\bigvee_{e\in E}\phi_e(x)$, for $E\subseteq \{0,1\}^k$,  is regular.

\begin{Que} \label{Q:gen-reg}
Is every property $\phi_E(x)$ of $L=\{U_1,\ldots,U_k\}$ regular?
\end{Que}

More generally:

\begin{Que} \label{Q:all-reg}
Is every property  of a finite relational language regular?
\end{Que}

\subsection{A necessary condition for regularity}

Next we shall  give a necessary condition in order for a property $\phi(x)$, (a) to have typicality degree $0$, and (b) to be regular.

Let $L$ be a language not necessarily relational. For any  $L$-property $\phi(x)$  and every $m\geq 1$ let us set
$$\phi^{(m)}:=(\exists x_1\cdots\exists x_m)\left((\bigwedge_{i\neq j}x_i\neq x_j)\wedge (\bigwedge_{i=1}^m\phi(x_i))\right).$$
$\phi^{(m)}$ is a sentence and says that $\phi(x)$ is satisfied by at least $m$ objects. So  $m<k$ implies  $\phi^{(k)}\rightarrow \phi^{(m)}$, and thus for every $m<k\leq n$,
\begin{equation} \label{E:smaller}
{\rm Mod}_n(\phi^{(k)})\subseteq  {\rm Mod}_n(\phi^{(m)}).
\end{equation}
As usual for every $n$ let $\lfloor\frac{n}{2}\rfloor$ be the  greatest integer $\leq n/2$. Then notice that  by the definition of $\textbf{S}_n(\phi:\mbox{typ})$,   for every  ${\cal M}\in \textbf{S}_n(L)$,
$${\cal M}\in \textbf{S}_n(\phi:\mbox{typ})\Leftrightarrow {\cal M}\models \phi^{(\lfloor\frac{n}{2}\rfloor+1)},$$
or
\begin{equation} \label{E:altern}
\textbf{S}_n(\phi:\mbox{typ})={\rm Mod}_n(\phi^{(\lfloor\frac{n}{2}\rfloor+1)}).
\end{equation}
On the other hand,  for every ${\cal M}$ and $k\leq |M|$, $|\phi({\cal M})|=k\Rightarrow {\cal M}\models \phi^{(k)}$, so
$${\cal M}\in \textbf{S}_{2n}(\phi:\mbox{ntr})\Rightarrow {\cal M}\models \phi^{(n)},$$
or
\begin{equation} \label{E:altern1}
\textbf{S}_{2n}(\phi:\mbox{ntr})\subseteq {\rm Mod}_{2n}(\phi^{(n)}).
\end{equation}

\begin{Lem} \label{L:key1}
Let $L$ be any  language and  $\phi(x)$ be an $L$-property. If there is $m\geq 1$ such that $\mu(\phi^{(m)})=0$, then:

(i)  $d(\phi:\mbox{typ})=0$, and

(ii) $\phi(x)$ is regular.
\end{Lem}

{\em Proof.}  Let $\mu(\phi^{(m)})=0$ for some fixed $m$. It means that
\begin{equation} \label{E:limit}
\mu_n(\phi^{(m)})=\frac{|{\rm Mod}_n(\phi^{(m)})|}{|\textbf{S}_n(L)|}\longrightarrow_n0.
\end{equation}

(i) For all $n\geq 2m$ we have  $m<\lfloor\frac{n}{2}\rfloor+1$, so by (\ref{E:smaller}) and (\ref{E:altern})
$$\textbf{S}_n(\phi:\mbox{typ})={\rm Mod}_n(\phi^{(\lfloor\frac{n}{2}\rfloor+1)})\subseteq  {\rm Mod}_n(\phi^{(m)}).$$
Consequently, for every $n\geq 2m$,
\begin{equation} \label{E:below}
\frac{|\textbf{S}_n(\phi:\mbox{typ})|}{|\textbf{S}_n(L)|}\leq \frac{|{\rm Mod}_n(\phi^{(m)})|}{|\textbf{S}_n(L)|}.
\end{equation}
Then (\ref{E:below}) combined with (\ref{E:limit}) yields
$$\frac{|\textbf{S}_n(\phi:\mbox{typ})|}{|\textbf{S}_n(L)|}\longrightarrow_n0,$$
i.e.  $d(\phi:\mbox{typ})=0$.

(ii) For every $n\geq m$, by (\ref{E:smaller}), (\ref{E:altern1}) and (\ref{E:limit}) we have
$$d_{2n}(\phi:\mbox{ntr})=\frac{|\textbf{S}_{2n}(\phi:\mbox{ntr})|}{|\textbf{S}_{2n}(L)|}\leq \frac{|{\rm Mod}_{2n}(\phi^{(n)})|}{|\textbf{S}_{2n}(L)|}\leq \frac{|{\rm Mod}_{2n}(\phi^{(m)})|}{|\textbf{S}_{2n}(L)|}\longrightarrow_n 0.$$
Thus $\phi(x)$ is regular.  \telos

\vskip 0.2in

\begin{Cor} \label{C:contra}
If $d(\phi:\mbox{typ})=1$, then $(\forall m\geq 1)(\mu(\phi^{(m)})=1)$.
\end{Cor}

{\em Proof.} Assume  $d(\phi:\mbox{typ})=1$. First note that if $L$ is relational, then the claim follows  immediately from \ref{L:key1} by the help of Theorem \ref{T:fagin} about the 0-1 law for  sentences of a relational $L$. For $d(\phi:\mbox{typ})=1$ implies $d(\phi:\mbox{typ})\neq 0$, so by \ref{L:key1} $(\forall m\geq 1)(\mu(\phi^{(m)})\neq 0)$ is true, and hence by \ref{T:fagin} $(\forall  m\geq 1)(\mu(\phi^{(m)})=1)$.

However the claim can be shown without appeal to Theorem \ref{T:fagin}, by a direct  argument similar to that of \ref{L:key1}. Namely, assuming $d(\phi:\mbox{typ})=1$ we have
\begin{equation} \label{E:converselimit}
\frac{|\textbf{S}_n(\phi:\mbox{typ})|}{|\textbf{S}_n(L)|}\longrightarrow_n1,
\end{equation}
and since by (\ref{E:below}) above,  we have that for every $n\geq 2m$,
$$\frac{|\textbf{S}_n(\phi:\mbox{typ})|}{|\textbf{S}_n(L)|}\leq \frac{|{\rm Mod}_n(\phi^{(m)})|}{|\textbf{S}_n(L)|}, $$
it follows that for all $m\geq 1$,
$$\frac{|{\rm Mod}_n(\phi^{(m)})|}{|\textbf{S}_n(L)|}\longrightarrow_n1,$$
i.e., $(\forall m\geq 1)(\mu(\phi^{(m)})=1)$.  \telos

\section{The 0-1 law for  typicality degrees and its failure for languages with unary properties}

Given a language $L$, the 0-1 law for typicality degrees of properties of $L$ can be defined by complete analogy with  the corresponding law for sentences described in Section 2.1, as follows.

\begin{Def} \label{D:degrees01}
{\em Let $L$ be a finite language. We say that} the 0-1 law holds for typicality degrees of properties of $L$, {\em if for every $L$-property $\phi(x)$, either $d(\phi:\mbox{typ})=0$ or $d(\phi:\mbox{typ})=1$.}
\end{Def}
We show  in this  Section  that the 0-1 law fails  for the   language  $L=\{U_1,\ldots,U_k\}$ with $k$ unary predicates. Recall from the previous section the definition of a basic  property $\phi_{\bar{s},e}(x)$ of $L$, for a $p$-subsequence $\bar{s}$ of $\langle 1,\ldots,k\rangle$ and an $e\in\{0,1\}^k$.

\begin{Thm} \label{T:counter-many}
Let $L=\{U_1,\ldots,U_k\}$ with $k\geq 1$, and $\phi_{\bar{s},e}(x)$ be a basic  property  of  $L$.

(i) If $\bar{s}$ is a $1$-sequence, then $d(\phi_{\bar{s},e}:\mbox{typ})=1/2$.

(ii) If $\bar{s}$ is a $p$-sequence for $p\geq 2$, then $d(\phi_{\bar{s},e}:\mbox{typ})=0$.
\end{Thm}

The subcase for $p=2$ of this Theorem requires special treatment, and it is more convenient to consider it separately. So  we shall split Theorem  \ref{T:counter-many} into Lemmas  \ref{L:case1} and \ref{L:case2} below. The proof of the Theorem is an immediate consequence of these Lemmas.

\begin{Lem} \label{L:case1}
(i) If $\bar{s}$ is a $1$-sequence, then $d(\phi_{\bar{s},e}:\mbox{typ})=1/2$.

(ii) If $\bar{s}$ is a $p$-sequence for $p\geq 3$, then $d(\phi_{\bar{s},e}:\mbox{typ})=0$.
\end{Lem}

{\em Proof.} Let  ${\cal M}=\langle M,W_1,\ldots,W_k\rangle$ be an $L$-structure with $|M|=n$, and fix a property $\phi_{\bar{s},e}(x)$ as above.  Recall that given $\phi_{\bar{s},e}({\cal M})=W^{e(i_1)}_{i_1}\cap \cdots\cap W^{e(i_p)}_{i_p}$. The proof has many similarities with the proof of Theorem \ref{T:regularity}. Given any set $A\subseteq M$ let us define $Z(A)$, $Z_1(A)$, $Z_2(A)$ exactly as in the aforementioned proof. Then, as before,   $|Z(A)|=|Z_1(A)|\cdot|Z_2(A)|$.

Setting $Z(m)=\bigcup\{Z(A):|A|=m\}$, we have $|Z(m)|=\binom{n}{m}\cdot|Z(A)|$, for any $A$ with $|A|=m$,  and also
\begin{equation} \label{E:sumup}
|\textbf{S}_n(\phi_{\bar{s},e}:\mbox{typ})|=\sum_{m>n/2}|Z(m)|.
\end{equation}

\vskip 0.1in

(i) Let $\bar{s}$ be a $1$-sequence. Without loss of generality we assume that $\phi_{\bar{s},e}(x)=U_1(x)$, so $\phi_{\bar{s},e}({\cal M})=W_1$.
Arguing exactly as in the proof of \ref{T:regularity}, we see that $|Z_1(A)|=1$, and
$$|Z_2(A)|= (2^n-1)_{k-1}\sim 2^{(k-1)n}.$$
Thus also $|Z(A)|=|Z_2(A)|\sim 2^{(k-1)n}$, and therefore
$$\textbf{S}_n(\phi_{\bar{s},e}(x):\mbox{typ})\sim \sum_{m>n/2}\binom{n}{m}\cdot
2^{(k-1)n}\sim 2^{n-1}\cdot 2^{(k-1)n},$$
hence  by (\ref{E:denom}),
$$d_n(\phi_{\bar{s},e}:\mbox{typ})\sim \frac{2^{n-1}\cdot 2^{(k-1)n}}{2^{kn}}\sim \frac{2^{n-1}}{2^n}=\frac{1}{2}.$$
Thus  $d(\phi_{\bar{s},e}:\mbox{typ})=\frac{1}{2}$, according to Fact \ref{F:asymlimit}.

\vskip 0.1in

(ii) Let $\bar{s}$ be a $p$-sequence with $3\leq p\leq k$.  Arguing as in the corresponding part of the proof of \ref{T:regularity},  for each $A$ with $|A|=m$,  the $p$-tuples of  $Z_1(A)$ are as many as the $p$-tuples $\langle Y_{j_1},\ldots,Y_{j_p}\rangle$ of elements of ${\cal P}(M\backslash A)$ whose members are pairwise disjoint. The number of all such sequences is $(2^{n-m})_p$ and constitutes an upper bound of $|Z_1(A)|$, that is
$$|Z_1(A)|\leq (2^{n-m})_{p}.$$
So
$$|Z_1(A)|\leq (2^{n-m})_p, \ \mbox{while} \ |Z_2(A)|=(2^n-p)_{k-p}.$$
Therefore
$$|Z(A)|=|Z_1(A)|\cdot|Z_2(A)|\leq(2^n-p)_{k-p}\cdot (2^{n-m})_p,$$ and
$$|Z(m)|\leq\binom{n}{m}\cdot (2^n-p)_{k-p}\cdot (2^{n-m})_p.$$
So
$$\textbf{S}_n(\phi_{\bar{s},e}:\mbox{typ})\leq \sum_{m>n/2}\binom{n}{m}[(2^n-p)_{k-p}\cdot  (2^{n-m})_p]=$$
$$=(2^n-p)_{k-p}\cdot \sum_{m>n/2}\binom{n}{m}(2^{n-m})_p \ \sim \ 2^{(k-p)n}\cdot\sum_{m>n/2}\binom{n}{m}(2^{n-m})_p.$$
Since $|\textbf{S}_n(L)|=(2^n)_k\ \sim \ 2^{kn}$, it follows from the last relation that
$$d_n(\phi_{\bar{s},e}:\mbox{typ})\lesssim  \frac{1}{2^{pn}}\cdot\sum_{m>n/2}\binom{n}{m}(2^{n-m})_p.$$
Setting $n-m=i$, this is written
\begin{equation} \label{E:exact?}
d_n(\phi_{\bar{s},e}:\mbox{typ})\lesssim  \frac{1}{2^{pn}}\cdot\sum_{i<n/2}\binom{n}{i}(2^i)_p.
\end{equation}
[Note that for $p>2^i$, i.e., for $i<\log_2p$, $(2^i)_p=0$, so the above sum is equal to $\sum_{\log_2p\leq i<n/2}\binom{n}{i}(2^i)_p$, however for notational simplicity we let $i$ range over all $i<n/2$.]

Now for every $n$, $k$, besides the relation $(n)_k\sim n^k$, the relation $(n)_k\leq n^k$ holds too. In particular  $(2^i)_p\leq 2^{pi}$, and also  for $i<n/2$, $2^{pi}<2^{\frac{pn}{2}}$, so (\ref{E:exact?}) implies
$$d_n(\phi_{\bar{s},e}:\mbox{typ})\lesssim  \frac{2^{\frac{pn}{2}}}{2^{pn}}\cdot\sum_{i<n/2}\binom{n}{i}=\frac{1}{2^{\frac{pn}{2}}}
\sum_{i<n/2}\binom{n}{i}.$$
Now by Fact \ref{F:semi-sum},  $\sum_{i<n/2}\binom{n}{i}\ \sim \ 2^{n-1}$, therefore
\begin{equation} \label{E:final}
d_n(\phi_{\bar{s},e}:\mbox{typ})\lesssim \frac{2^{n-1}}{2^{\frac{pn}{2}}}.
\end{equation}
Since  $p\geq 3$, we have
$$\frac{2^{n-1}}{2^{\frac{pn}{2}}}\leq \frac{2^{n-1}}{2^{\frac{3n}{2}}}=\frac{1}{2^{\frac{n}{2}+1}}\longrightarrow_n0,$$
so, by Fact \ref{F:asymlimit},  (\ref{E:final}) yields  $d_n(\phi_{\bar{s},e}:\mbox{typ})\longrightarrow_n0$. This proves clause (ii) and completes the proof of the Lemma.  \telos

\vskip 0.2in

In the proof of the next Lemma we shall make use of  Stirling numbers of the second kind.  Recall that   they are denoted ${n \brace k}$, or $S(n,k)$,  where,   for $1\leq k\leq n$,  ${n \brace k}$ counts the number  of partitions of $\{1,\ldots,n\}$ into $k$ nonempty parts. The  explicit formula for ${n \brace k}$ is  (see \cite[p. 231]{Hetal08} or \cite{Stirling}):
\begin{equation} \label{E:explicit}
{n \brace k}=\frac{1}{k!}\sum_{i=0}^k(-1)^{k-i}\binom{k}{i}i^n.
\end{equation}

\begin{Lem} \label{L:case2}
If $\bar{s}$ is a $2$-sequence  then $d(\phi_{\bar{s},e}:\mbox{typ})=0$.
\end{Lem}

{\em Proof.} Fix a $2$-subsequence  $\bar{s}=\langle i_1,i_2\rangle$ of $\langle  1,\ldots,k\rangle$ and an $e\in \{0,1\}^k$. For any set $A\subseteq M$,
$$Z_1(A)=\{\langle W_{i_1}^{e(i_1)},W_{i_2}^{e(i_2)}\rangle: W_{i_1}^{e(i_1)}\cap W_{i_2}^{e(i_2)}=A\},$$
so $W_{i_1}^{e(i_1)}=A\cup Y_1$ and  $W_{i_2}^{e(i_2)}=A\cup Y_2$, where $Y_1,Y_2\subseteq M\backslash A$ and $Y_1\cap Y_2=\emptyset$. It follows that $|Z_1(A)|=|P(A)|$, where
$$P(A)=\{\langle Y_1,Y_2\rangle: Y_1,Y_2\subseteq M\backslash A \ \wedge \ Y_1\cap Y_2=\emptyset\}.$$
The pairs $\langle Y_1,Y_2\rangle$ are of two kinds: those for which $Y_1,Y_2$ form a partition of $M\backslash A$, i.e.,  $Y_1\cup Y_2=M\backslash A$, and those for which $Y_1\cup Y_2\neq M\backslash A$. That is, $P(A)=P_1(A)\cup P_2(A)$, where
$$P_1(A)=\{\langle Y_1,Y_2\rangle: Y_1\cap Y_2=\emptyset \ \wedge \ Y_1\cup Y_2= M\backslash A\},$$
$$P_2(A)=\{\langle Y_1,Y_2\rangle: Y_1\cap Y_2=\emptyset \ \wedge \ Y_1\cup Y_2\neq  M\backslash A\},$$
and
\begin{equation} \label{E:part}
|Z_1(A)|=|P(A)|=|P_1(A)|+|P_2(A)|.
\end{equation}
Now $|P_1(A)|$ and $|P_2(A)|$ can be easily calculated in terms of the  $2$-partitions and  $3$-partitions of  $M\backslash A$, respectively. Let $\Pi(M\backslash A,2)$, $\Pi(M\backslash A,3)$ denote  the sets of partitions of $M\backslash A$ into 2 and 3 nonempty parts, respectively. Then $|\Pi(M\backslash A,2)|={n-m \brace 2}$ and $|\Pi(M\backslash A,3)|={n-m \brace 3}$. Now a member of $\Pi(M\backslash A,2)$ is a 2-element subset of $M\backslash A$, while  $P_1(A)$ consists of {\em ordered pairs} of such subsets, therefore
$$|P_1(A)|=2\cdot |\Pi(M\backslash A,2)|=2\cdot{n-m \brace 2}.$$
Analogously every member of $\Pi(M\backslash A,3)$ is a 3-element subset of $M\backslash A$, and each such subset provides $(3)_2$ pairs that belong to  $P_2(A)$ so
$$|P_2(A)|=(3)_2\cdot |\Pi(M\backslash A,3)|=6\cdot{n-m \brace 3}.$$
For every $n\geq 2$, it is easy to see (without appealing to (\ref{E:explicit})) that ${n \brace 2}=2^{n-1}-1$, so  for every $A$ with  $|A|=m\leq n-2$,
$$|P_1(A)|=2\cdot (2^{n-m-1}-1)=2^{n-m}-2.$$
If $P_1(m)=\bigcup\{P_1(A):|A|=m\}$, then for $m\leq n-2$,
\begin{equation} \label{E:part1}
|P_1(m)|=\binom{n}{m}\cdot (2^{n-m}-2).
\end{equation}
On the other hand, for every $n\geq 3$,   ${n \brace 3}$ is calculated by the help of    formula (\ref{E:explicit}) which yields:
$${n \brace 3}=\frac{1}{6}[3-3\cdot 2^n+3^n].$$
Therefore,  for $|A|=m\leq n-3$,
$$|P_2(A)|=3-3\cdot 2^{n-m}+3^{n-m}.$$
So setting as before $P_2(m)=\bigcup\{P_2(A):|A|=m\}$, we have for $m\leq n-3$,
\begin{equation} \label{E:part2}
|P_2(m)|=\binom{n}{m}\cdot (3-3\cdot 2^{n-m}+3^{n-m}).
\end{equation}
Finally by (\ref{E:part}), (\ref{E:part1}) and (\ref{E:part2}) above we obtain, for $|A|=m\leq n-3$,
$$|Z_1(m)|=|P_1(m)|+|P_2(m)|=\binom{n}{m}(3^{n-m}-2^{n-m+1}+1).$$
Recall also from above that for $p=2$,
$$|Z_2(m)|\sim \binom{n}{m} 2^{(k-2)n},$$
so for $m\leq n-3$,
\begin{equation} \label{E:omit}
|Z(m)|=|Z_1(m)|\cdot|Z_2(m)|\sim \binom{n}{m}\cdot[2^{(k-2)n} \cdot (3^{n-m}-2^{n-m+1}+1)].
\end{equation}
Now it is easy to see that
$$\textbf{S}_n(\phi_{\bar{s},e}:\mbox{typ})=\sum_{n\geq m>n/2}|Z(m)|\sim  \sum_{n-3\geq m>n/2}|Z(m)|,$$
so by (\ref{E:omit}),
$$\textbf{S}_n(\phi_{\bar{s},e}:\mbox{typ})\sim 2^{(k-2)n}\cdot\sum_{n-3\geq m>n/2}
\binom{n}{m}\cdot (3^{n-m}-2^{n-m+1}+1),$$
and, given that $\textbf{S}_n(L)=2^{kn}$,
$$d_n(\phi_{\bar{s},e}:\mbox{typ})\sim \frac{1}{2^{2n}}\cdot\sum_{n-3\geq m>n/2}
\binom{n}{m}\cdot (3^{n-m}-2^{n-m+1}+1).$$
Setting $n-m=i$, this is written
$$d_n(\phi_{\bar{s},e}:\mbox{typ})\sim \frac{1}{2^{2n}}\cdot\sum_{3\leq i<n/2}
\binom{n}{i}\cdot (3^i-2^{i+1}+1)\leq\frac{1}{2^{2n}}\cdot\sum_{3\leq i<n/2}
\binom{n}{i}\cdot (3^{\frac{n}{2}}-2^{\frac{n}{2}+1}+1)=$$
$$=\frac{3^{\frac{n}{2}}-2^{\frac{n}{2}+1}+1}{2^{2n}}\cdot \sum_{3\leq i<n/2}
\binom{n}{i}\leq\frac{3^{\frac{n}{2}}-2^{\frac{n}{2}+1}+1}{2^{2n}}\cdot \sum_{0\leq i<n/2}\binom{n}{i}.$$
By Fact \ref{F:semi-sum}, the last quantity is
$$\sim \frac{3^{\frac{n}{2}}-2^{\frac{n}{2}+1}+1}{2^{2n}}\cdot 2^{n-1}\sim \frac{3^{\frac{n}{2}}-2^{\frac{n}{2}+1}+1}{2^{n+1}}\sim \frac{1}{2}\left(\frac{3^{\frac{n}{2}}}{4^{\frac{n}{2}}}-
\frac{1}{2^{\frac{n}{2}-1}}+\frac{1}{2^n}\right).$$
So finally,
$$d_n(\phi_{\bar{s},e}:\mbox{typ})\lesssim  \frac{1}{2}\left(\frac{3^{\frac{n}{2}}}{4^{\frac{n}{2}}}-
\frac{1}{2^{\frac{n}{2}-1}}+\frac{1}{2^n}\right)
\longrightarrow_n0.$$
This completes the proof. \telos

\vskip 0.2in

An immediate  consequence of clause (i) of Theorem \ref{T:counter-many} is the following.

\begin{Cor} \label{C:0-1fails}
The 0-1 law  for  typicality degrees of properties of a relational language fails in general.
\end{Cor}

Let us consider at this point  another question related to   property ``$U(x)$'', namely the question about the probability of the sentences $U(x)^{(m)}$, for $m\geq 1$. We can see that, in contrast to the typicality degree $1/2$ of $U(x)$, the truth probability of the sentences $U(x)^{(m)}$ is $1$.

\begin{Prop} \label{P:probability}
For every  $m\geq 1$, $\mu(U(x)^{(m)})=1$.
\end{Prop}

{\em Proof.} The sentence $U(x)^{(m)}$ says that ``$U(x)$ is satisfied by at least $m$ elements''. Therefore ${\rm Mod}_n(U(x)^{(m)})=\{A\subseteq\{1,\ldots,n\}:|A|\geq m\}$, or  $|{\rm Mod}_n(U(x)^{(m)})|=2^n-\sum_{0\leq i<m}\binom{n}{i}$, and therefore
$$\mu_n(U(x)^{(m)})=\frac{2^n-\sum_{0\leq i<m}\binom{n}{i}}{2^n}=1-\frac{\sum_{0\leq i<m}\binom{n}{i}}{2^n}.$$
The nominator $\sum_{0\leq i<m}\binom{n}{i}$ of the fraction on the right-hand side is a polynomial in $n$ of degree $m-1$, so its quotient by the exponential $2^n$ goes to $0$ an $n$ grows. Thus $\mu(U(x)^{(m)})=\lim_{n\rightarrow\infty}\mu_n(U(x)^{(m)})=1$. \telos

\begin{Que} \label{Q:zero-one}
Does the 0-1 law about typicality degrees hold for every property of a relational language  without unary predicates?
\end{Que}

We close this section with a remark about Question \ref{Q:zero-one}. Some people believe that  the answer to this  question must be  positive on the basis that  the  rather general method of ``extension axioms'' that was  used by R. Fagin in \cite{Fa76} to prove the 0-1 law for truth degrees of {\em sentences}  of a relational language, could be also applied somehow to the case of typicality degrees of properties.\footnote{In contrast the method of Glebskii {\em et al.} in \cite{GETAL69} seems to be rather  ad hoc.} The problem however is that this specific method works for sentences of {\em every} relational language, including  those that contain  unary predicates, while, as we saw above,  in the case of typicality  degrees the method should not work when unary predicates are included. I don't know how this gap could  be bridged and if Fagin's method is  actually applicable to the present case.

\section{Some results about regularity and  degrees of properties of the  language $\textbf{L=\{F\}}$}

The failure of 0-1 law for typicality degrees of properties of languages with unary predicates is a divergence from the behavior of truth probabilities of sentences.  In this Section we consider some properties of the language $L=\{F\}$ where $F$ is a unary function symbol. We saw in  Example \ref{E:example1} that $\mu(\forall x (F(x)\neq x))=e^{-1}$, which means that the 0-1 law does not hold for sentences of $L$. The question is whether the 0-1 law fails also for typicality degrees of properties of this language. In this section we consider two properties: 1) $\phi(x):=(F(x)\neq x)$ and 2) $\psi(x):=\exists y(F(y)=x)$. It is shown that both are regular, the degree of $\phi(x)$ is 1, while the degree of $\psi(x)$ is not known, although we give evidence that it is 1 too.

\subsection{The property $\phi(\textbf{x):=(F(x)}\neq \textbf{x})$ and its negation}

\begin{Prop} \label{P:counter}
Let $L=\{F\}$ and let $\phi(x):=(F(x)\neq x)$. Then $$d(\phi:\mbox{typ})=1.$$
\end{Prop}

{\em Proof.} Let ${\cal M}=\langle M,f\rangle$ with $|M|=n$ and $A\subseteq M$ such that $A=\phi({\cal M})$. Then $a\in A\Leftrightarrow f(a)\neq a$ and $a\notin A\Leftrightarrow f(a)=a$. Thus if $G(A)=\{f\in M^M:\phi({\cal M})=A\}$, every $f\in G(A)$ can be identified with  the pair $(f\restr A, id_{(M\backslash A)})$, or since $id_{(M\backslash A)}$ is unique, with $f\restr (M\backslash A)$. As we argued in Example \ref{E:example1},  if $|A|=m$ then $|G(A)|=(n-1)^m$. Let us set  $G(m)=\bigcup\{G(A):A\subseteq M \ \& \ |A|=m\}$.  Then, since there are  $\binom{n}{m}$ sets of cardinality $m$, $|G(m)|=\binom{n}{m}(n-1)^m$. Also by definition,
$$\textbf{S}_n(\phi:\mbox{typ})=\bigcup\{G(m):m>n/2\},$$
therefore
\begin{equation} \label{E:atlast}
|\textbf{S}_n(\phi:\mbox{typ})|=\sum_{m>n/2}|G(m)|=\sum_{m>n/2}\binom{n}{m}(n-1)^m.
\end{equation}
It is more convenient to set $m=n-k$ and write this sum in the form:
$$|\textbf{S}_n(\phi:\mbox{typ})|=\sum_{k<n/2}\binom{n}{n-k}(n-1)^{n-k}=
\sum_{k<n/2}\binom{n}{k}(n-1)^{n-k},$$
or
\begin{equation} \label{E:coef}
|\textbf{S}_n(\phi:\mbox{typ})|=\sum_{k<n/2}a_k(n),
\end{equation}
where $a_k(n)=\binom{n}{k}(n-1)^{n-k}$. Then
$$d_n(\phi:\mbox{typ})=\frac{|\textbf{S}_n(\phi:\mbox{typ})|}{n^{n}}=
\sum_{k<n/2}\frac{a_k(n)}{n^{n}}.$$
Setting for simplicity,
$a_n=\sum_{k<n/2}\frac{a_k(n)}{n^{n}}$, we have to compute the limit
\begin{equation} \label{E:first}
d(\phi:\mbox{typ})=\lim_{n\rightarrow\infty}a_n.
\end{equation}
Now $\frac{a_k(n)}{n^n}=\frac{1}{n^n}\cdot \binom{n}{k}(n-1)^{n-k}$,
or, multiplying and dividing the right-hand side by  $(n-1)^k$, this is written
$$\frac{a_k(n)}{n^{n}}=\frac{(n-1)^n}{n^n}\cdot \binom{n}{k}\frac{1}{(n-1)^k}.$$
Therefore
$$a_n=\sum_{k<n/2}\left(\frac{n-1}{n}\right)^n\cdot \binom{n}{k}\frac{1}{(n-1)^k}=
\left(\frac{n-1}{n}\right)^n\cdot\sum_{k<n/2} \binom{n}{k}\frac{1}{(n-1)^k}=b_n\cdot c_n,$$
where  $b_n=\left(\frac{n-1}{n}\right)^n$ and $c_n=\sum_{k<n/2} \binom{n}{k}\frac{1}{(n-1)^k}$, and
\begin{equation} \label{E:second}
\lim_{n\rightarrow\infty}a_n=(\lim_{n\rightarrow\infty}b_n)\cdot(\lim_{n\rightarrow\infty}c_n).
\end{equation}
Now $\lim_{n\rightarrow\infty}b_n=e^{-1}$. So it remains to compute
\begin{equation} \label{E:third}
\lim_{n\rightarrow\infty}c_n=\lim_{n\rightarrow\infty}\sum_{k<n/2}
\binom{n}{k}\frac{1}{(n-1)^k}.
\end{equation}
We shall show that $\lim_{n\rightarrow\infty}c_n=e$. It suffices to show that $$\lim_{n\rightarrow\infty}c_{2n}=\lim_{n\rightarrow\infty}c_{2n+1}=e.$$

\vskip 0.1in

(a) \underline{Proof of $\lim_{n\rightarrow\infty}c_{2n}=e$.}  We have
$$c_{2n}=\sum_{k<n}\binom{2n}{k}\frac{1}{(2n-1)^k}=
\sum_{k=0}^{n-1}\binom{2n}{k}\frac{1}{(2n-1)^k}.$$
However it is easy to see that  $\lim_{n\rightarrow\infty}c_{2n}=\lim_{n\rightarrow\infty}c'_{2n}$, where
$$c'_{2n}=\sum_{k=0}^n\binom{2n}{k}\frac{1}{(2n-1)^k}.$$
This is because  $c'_{2n}-c_{2n}=\binom{2n}{n}\frac{1}{(2n-1)^n}$, which goes to $0$ when $n$ goes to infinity, as follows easily from Fact \ref{F:central}.  Therefore it suffices to show that
$$\lim_{n\rightarrow\infty}c'_{2n}=e.$$
To show that, we  compare each term $c'_{2n}$ with the term $$g_n=\sum_{k=0}^n\binom{n}{k}\frac{1}{(n-1)^k}=\left(1+\frac{1}{n-1}\right)^n=
\left(\frac{n}{n-1}\right)^n,$$ for which it is well-known that   $\lim_{n\rightarrow\infty}g_n=e$. Specifically we compare each summand $A_k(n)=\binom{2n}{k}\frac{1}{(2n-1)^k}$ of $c_{2n}'$, for $k\leq n$,  with the corresponding summand $B_k(n)=\binom{n}{k}\frac{1}{(n-1)^k}$ of $g_n$. Then we have
$$A_k(n)=\frac{1}{k!}\frac{(2n-k+1)(2n-k+2)\cdots 2n}{(2n-1)^k}=\frac{1}{k!}\frac{P(2n)}{Q(2n)},$$
while
$$B_k(n)=\frac{1}{k!}\frac{(n-k+1)(n-k+2)\cdots n}{(n-1)^k}=\frac{1}{k!}\frac{P(n)}{Q(n)},$$
where $P(x)$ and $Q(x)$ are polynomials of degree $k$. If $\alpha_k$ and $\beta_k$ are    the leading coefficients of $P(x)$ and $Q(x)$, respectively, then $$\lim_{n\rightarrow\infty}\frac{P(2n)}{Q(2n)}=
\lim_{n\rightarrow\infty}\frac{P(n)}{Q(n)}=\frac{\alpha_k}{\beta_k},$$ therefore $\lim_{n\rightarrow\infty}(A_k(n)-B_k(n))=0$. Besides,
$$c'_{2n}-g_n=\sum_{k=0}^nA_k(n)-\sum_{k=0}^nB_k(n)=\sum_{k=0}^n(A_k(n)-B_k(n)),$$
so
$$\lim_{n\rightarrow\infty}(c'_{2n}-g_n)=
\sum_{k=0}^n\lim_{n\rightarrow\infty}(A_k(n)-B_k(n))=0.$$
Therefore $\lim_{n\rightarrow\infty}c'_{2n}=\lim_{n\rightarrow\infty}g_n=e$, as required.
\vskip 0.1in

(b) \underline{Proof of $\lim_{n\rightarrow\infty}c_{2n+1}=e$.} This is almost the same as in (a). First notice  that $k<(2n+1)/2 \Leftrightarrow k\leq n$, so
$$c_{2n+1}=\sum_{k<(2n+1)/2}\binom{2n+1}{k}\frac{1}{(2n)^k}=
\sum_{k=0}^n\binom{2n+1}{k}\frac{1}{(2n)^k}.$$
Then we set for every $k\leq n$,   $$A_k'(n)=\binom{2n+1}{k}\frac{1}{(2n)^k}=\frac{1}{k!}\frac{P(2n+1)}{Q(2n+1)}$$ and compare it again with
$$B_k(n)=\frac{1}{k!}\frac{P(n)}{Q(n)}$$
of the preceding case. As before we have
$$\lim_{n\rightarrow\infty}\frac{P(2n+1)}{Q(2n+1)}=
\lim_{n\rightarrow\infty}\frac{P(n)}{Q(n)}=\frac{\alpha_k}{\beta_k},$$ therefore $\lim_{n\rightarrow\infty}(A'_k(n)-B_k(n))=0$, and hence
$$\lim_{n\rightarrow\infty}(c_{2n+1}-g_n)=
\sum_{k=0}^n\lim_{n\rightarrow\infty}(A'_k(n)-B_k(n))=0.$$
So $\lim_{n\rightarrow\infty}c_{2n+1}=\lim_{n\rightarrow\infty}g_n=e$. This completes the proof. \telos

\begin{Prop} \label{P:regularf}
The property $\phi(x):=(F(x)\neq x)$ is regular.
\end{Prop}

{\em Proof.} Let ${\cal M}=\langle M,f\rangle$ be an $L$-structure with $|M|=2n$. Then ${\cal M}\in \textbf{S}_{2n}(\phi:\mbox{ntr})$ if and only if  $|\phi({\cal M})|=n$. Let $A=\phi({\cal M})$. Then $A=\{a\in M:f(a)\neq a\}$ and hence $a\notin A\Leftrightarrow f(a)=a$, that is, $f\restr(M\backslash A)=id_{M\backslash A}$. Since $id_{M\backslash A}$ is unique, it follows that if $G(A)$ is the set of $f\in M^M$ such that $\phi({\cal M})=A$, then
$$|G(A)|=|\{f\restr A: \phi({\cal M})=A\}|=|\{g\in M^A: \forall x (g(x)\neq x)\}|.$$ Since $|M|=2n$ and $|A|=n$ we have,  as argued in Example \ref{E:example1},    $|G(A)|=(2n-1)^n$. Since there are $\binom{2n}{n}$ such sets $A$, it follows that   $$|\textbf{S}_{2n}(\phi:\mbox{ntr})|=|\bigcup\{G(A):|A|=n\}|=
\binom{2n}{n}\cdot(2n-1)^n.$$
On the other hand $|\textbf{S}_{2n}(L)|=(2n)^{2n}$, hence
$$\frac{|\textbf{S}_{2n}(\phi:\mbox{ntr})|}{|\textbf{S}_{2n}(L)|}=
\frac{\binom{2n}{n}\cdot(2n-1)^n}{(2n)^{2n}}\leq \frac{\binom{2n}{n}\cdot(2n)^n}{(2n)^{2n}}=\frac{\binom{2n}{n}}{(2n)^n}.$$
By Fact \ref{F:central},  $\binom{2n}{n}\leq\frac{4^n}{\sqrt{\pi n}}$, so the preceding inequality implies
$$\frac{|\textbf{S}_{2n}(\phi:\mbox{ntr})|}{|\textbf{S}_{2n}(L)|}\leq \frac{4^n}{(2n)^n\cdot\sqrt{\pi n}}=
\frac{2^n}{n^n\cdot \sqrt{\pi n}}\longrightarrow_n 0.$$
\telos

\begin{Cor} \label{C:dual}
$d(F(x)=x:\mbox{typ})=0$.
\end{Cor}

{\em Proof.} By Proposition \ref{P:regularf}, $F(x)= x$ is regular, while by Proposition \ref{P:counter} $d(F(x)\neq x:\mbox{typ})=1$.  Therefore  in view of Fact \ref{F:sum} (iii),   $d(F(x)=x:\mbox{typ})=1-d(F(x)\neq x:\mbox{typ})=0$. \telos

\vskip 0.2in

Next Proposition shows  that the criterion of Lemma \ref{L:key1} (ii), which allows one to deduce the regularity of  property $\phi(x)$ whenever  for some $m\geq 1$, $\mu(\phi^{(m)})=0$,  cannot be used for properties  $F(x)\neq x$ and $F(x)=x$.

\begin{Prop} \label{P:applykey1}
For every $m\geq 1$,

(i) $\mu((F(x)\neq x)^{(m)})=1$.

(ii) $\frac{e^{-1}}{m!}\leq((F(x)=x)^{(m)})\leq \frac{1}{m!}$. In particular, $\mu((F(x)=x)^{(1)})=1-e^{-1}$.
\end{Prop}

{\em Proof.} (i) By definition
$(F(x)\neq x)^{(m)}$ holds in $\langle M,f\rangle$ if and only if  there are distinct $x_1,\ldots,x_m\in M$ such that $f(x_i)\neq x_i$ for every $i=1,\ldots,m$. Equivalently, if and only if  for every $A\subseteq M$ such that $f\restr A=id$, $|A|\leq |M|-m$. Fixing $M$ with $|M|=n$, for a given $A\subseteq M$ the totality of $f$ such that $f\restr A=id$ are $n^{n-|A|}$, while when $A$ ranges over all subsets with $|A|\leq n-m$, the totality of $f$ for which $\langle M,f\rangle$  satisfies $(F(x)\neq x)^{(m)}$ has cardinality
$$|{\rm Mod}_n((F(x)\neq x)^{(m)})|=\sum_{|A|\leq n-m}n^{n-|A|}=\sum_{m\leq i\leq n}n^i.$$
Therefore $\mu_n((F(x)\neq x)^{(m)})=\frac{\sum_{{m\leq i\leq n}}n^i}{n^n}\longrightarrow_n1.$

(ii) Let us show first the second claim,  that  $\mu((F(x)=x)^{(1)})=1-e^{-1}$. This follows from Example \ref{E:example1},  and the fact that for every sentence $\phi$, $\mu(\neg\phi)=1-\mu(\phi)$. So
$$\mu((F(x)=x)^{(1)})=\mu(\exists x(F(x)=x))=\mu(\neg (\forall x)(F(x)\neq x))=$$
$$=1-\mu(\forall x(F(x)\neq x))=1-\frac{1}{e}=\frac{e-1}{e}.$$

Now consider the sentence $(F(x)=x)^{(m)}$ for $m\geq 1$.  $(F(x)=x)^{(m)}$ holds in $\langle M,f\rangle$ if and only if  there are distinct $x_1,\ldots,x_m\in M$ such that $f(x_i)=x_i$ for every $i=1,\ldots,m$, i.e., if there is $A\subseteq M$ with $|A|=m$ such that $f\restr A=id$. If $|M|=n$, for each $A$ with $|A|=m$ there are $n^{n-m}$ functions such that $f\restr A=id$. Since there are $\binom{n}{m}$ such sets $A$ with $|A|=m$, the totality of functions of this kind is at most $\binom{n}{m}n^{n-m}$ (because  this  totality possibly contains repetitions). So $$|{\rm Mod}_n((F(x)= x)^{(m)})|\leq\binom{n}{m}n^{n-m}.$$
On the other hand, for a fixed $A$ with $|A|=m$, let $X_A$ be the collection of functions $f$ such that $f\restr A=id$ and $f(x)\neq x$ for every $x\in M\backslash A$. Since for every $f\in X_A$, $f(x)$ may take independently $n-1$ possible values on $M\backslash A$, it follows that $|X_A|=(n-1)^{n-m}$ and, moreover,  if $A\neq A'$ then $X_A\cap X_{A'}=\emptyset$. So
$$|{\rm Mod}_n((F(x)= x)^{(m)})|\geq\binom{n}{m}(n-1)^{n-m}.$$
Therefore
\begin{equation} \label{E:lb}
\frac{\binom{n}{m}(n-1)^{n-m}}{n^n}\leq \mu_n((F(x)=x)^{(m)})\leq \frac{\binom{n}{m}n^{n-m}}{n^n}.
\end{equation}
Denoting by $l_n$ the preceding lower bound in (\ref{E:lb}) we have $$l_n=\frac{\binom{n}{m}(n-1)^{n-m}}{n^n}=\frac{\binom{n}{m}}{(n-1)^m}\cdot \frac{(n-1)^n}{n^n}=$$
$$=\frac{1}{m!}\frac{(n-m+1)(n-m+2)\cdots (n-1)n}{(n-1)^m}\cdot \left(\frac{n-1}{n}\right)^n.$$
So
$$\lim_{n\rightarrow\infty} l_n=\frac{1}{m!}\lim_{n\rightarrow\infty}\frac{(n-m+1)(n-m+2)\cdots (n-1)n}{(n-1)^m}\cdot \lim_{n\rightarrow\infty}\left(\frac{n-1}{n}\right)^n.$$
Now $\lim_{n\rightarrow\infty}\frac{(n-m+1)(n-m+2)\cdots (n-1)n}{(n-1)^m}=1$, because both the nominator and the denominator are polynomials of degree $m$ with leading coefficients 1, while  $\lim_{n\rightarrow\infty} \left(\frac{n-1}{n}\right)^n=e^{-1}$. Therefore $\lim_{n\rightarrow\infty}l_n=\frac{e^{-1}}{m!}$.

Similarly, if $u_n=\frac{\binom{n}{m}n^{n-m}}{n^n}$ is the upper bound in (\ref{E:lb}) above, then
$$u_n=\frac{\binom{n}{m}}{n^m}=\frac{1}{m!}\cdot\frac{(n-m+1)(n-m+2)\cdots (n-1)n}{n^m},$$
and comparing as  before the polynomials of the fraction we find  $\lim_{n\rightarrow\infty}u_n=\frac{1}{m!}$. So finally
$$\frac{e^{-1}}{m!}\leq \mu((F(x)=x)^{(m)})=\lim_{n\rightarrow\infty}\mu_n((F(x)=x)^{(m)})\leq \frac{1}{m!}.$$
Notice that the value $1-e^{-1}=\frac{e-1}{e}$ of $\mu((F(x)=x)^{(1)})$ conforms with the general bounds given  above, since  $e^{-1}\leq \frac{e-1}{e}\leq 1$.  \telos

\vskip 0.2in

\begin{Cor} \label{C:noconverse}
The converse of Lemma \ref{L:key1} (i) is false for the language $L=\{F\}$. Namely, there is $\phi(x)$ of $L$ such that $d(\phi:\mbox{typ})=0$ while $\forall m\geq 1 \ \mu(\phi^{(m)})>0$.
\end{Cor}

{\em Proof.} Take $\phi(x):(F(x)=x)$. By Corollary \ref{C:dual},
$d(F(x)=x:\mbox{typ})=0$, while by Proposition  \ref{P:applykey1}, for every $m\geq 1$  $\mu((F(x)=x)^{(m)})\geq \frac{e^{-1}}{m!}$. \telos

\section{Degrees of some properties of graphs}
In this Section we examine the typicality degree of some natural first-order properties of finite undirected graphs. The language of graphs is $L=\{E\}$, where $E$ is the symbol of a binary symmetric and irreflexive relation. Thus an $L$-structure is a pair $G=\langle A,E\rangle$, where $A$ is the set of nodes of $G$ and $E$ is its  {\em adjacency}  relation.  $xEy$ means that  ``the nodes $x,y$ are adjacent'', i.e.,  connected  with an edge. By assumption $xEy\Leftrightarrow yEx$ and for all $x\in A$ $\neg (xEx)$.  So $E$  is a set of 2-element subsets of $A$, that is $E\subseteq [A]^2$, and $G\models xEy \Leftrightarrow \{x,y\}\in E$. So if $|A|=n$, then $|[A]^2|=\binom{n}{2}$, therefore  $|\textbf{S}_n(L)|=2^{\binom{n}{2}}$.

Let us consider  the following  examples of first-order properties of $L$:

\vskip 0.1in

(1) $\phi_{\rm none}(x)$: ``$x$ is an isolated node'' [($\forall y)\neg(xEy))$].

(2) $\phi_{\rm all}(x)$: ``$x$ is adjacent to every  node'' [$(\forall y)(xEy)$].

(3) $\phi_{\rm one}(x)$: ``$x$ is adjacent to exactly one  node'' [$(\exists!y)(xEy)$].

\begin{Prop} \label{P:application}
Properties $\phi_{\rm none}(x)$,  $\phi_{\rm all}(x)$ and $\phi_{\rm one}(x)$ have typicality degree $0$. That is,  $d(\phi_{\rm none}:\mbox{typ})=d(\phi_{\rm all}:\mbox{typ})=d(\phi_{\rm one}:\mbox{typ})=0$.
\end{Prop}

{\em Proof.} We shall prove the claim using Lemma \ref{L:key1} (i), namely it suffices to prove that $\mu(\phi_{\rm none}^{(1)})=\mu(\phi_{\rm all}^{(1)})=\mu(\phi_{\rm one}^{(1)})=0$.

(i) $\mu(\phi_{\rm none}^{(1)})=0$:  We have $\phi_{\rm none}^{(1)}:=(\exists x)(\forall y)\neg(xEy)$. We must show that
\begin{equation} \label{E:showlimit}
\frac{|{\rm Mod}_n(\phi_{\rm none}^{(1)})|}{|\textbf{S}_n(L)|}\longrightarrow_n 0.
\end{equation}
We saw above that $|\textbf{S}_n(L)|=2^{\binom{n}{2}}$, so we have to estimate also $|{\rm Mod}_n(\phi_{\rm none}^{(1)})|$. Let $G=\langle A,E\rangle$ be a graph with $n$ nodes satisfying $\phi_{\rm none}^{(1)}$. Then $G$  contains at {\em least one}  isolated node $a$ (see Figure 1).
\begin{center}
\setlength{\unitlength}{1mm}
\begin{picture}(30,35)(-15,-35)
\linethickness{1pt} \thinlines

\put(-40,0){\line(1,0){70}}
\put(-40,-35){\line(1,0){70}}
\put(-40,0){\line(0,-1){35}}
\put(30,0){\line(0,-1){35}}

\put(28,-2){\makebox(0,0)[c]{$G$}}
\put(-18,-17){\makebox(0,0)[c]{$a$}}
\put(-25,-30){\makebox(0,0)[c]{$A'$}}

\put(-15,-15){\circle*{2}}
\put(-15,-15){\circle{12}}

\end{picture}
\end{center}
\begin{center}
Figure 1
\end{center}

\vskip 0.1in

\noindent It follows that if $A'=A\backslash\{a\}$ and $G'=\langle A',E'\rangle$,  where $E'$ is the restriction of $E$ to $A'$, then $G'$ can be {\em any}  graph with $n-1$ nodes (including those with any number of  isolated nodes). Therefore there are $2^{\binom{n-1}{2}}$ different graphs on $A$ in which $a$ is an isolated node. Letting $a$ range over every node of $A$, it follows  that the total number of graphs in $\textbf{S}_n(L)$ that satisfy  $\phi_{\rm none}^{(1)}$, i.e., $|{\rm Mod}_n(\phi_{\rm none}^{(1)})|$, is at most $n 2^{\binom{n-1}{2}}$. Thus
$$\frac{|{\rm Mod}_n(\phi_{\rm none}^{(1)})|}{|\textbf{S}_n(L)|}\leq\frac{n2^{\binom{n-1}{2}}}
{2^{\binom{n}{2}}}=\frac{n2^{\frac{(n-2)(n-1)}{2}}}{2^{\frac{n(n-1)}{2}}}=
\frac{n}{2^{n-1}}\longrightarrow_n 0,$$
so (\ref{E:showlimit}) is true.

\vskip 0.1in

(ii)  $\mu(\phi_{\rm all}^{(1)})=0$:  We have $\phi_{\rm all}^{(1)}:=(\exists x)(\forall y)(xEy)$ and we must  show that
\begin{equation} \label{E:showlimit2}
\frac{|{\rm Mod}_n(\phi_{\rm all}^{(1)})|}{|\textbf{S}_n(L)|}\longrightarrow_n 0.
\end{equation}
The argument is quite  similar to that of the previous case. If $G=\langle A,E\rangle$ is a graph with $n$ nodes  satisfying $\phi_{\rm all}^{(1)}$, there is $a\in A$ which is adjacent to all other nodes of $G$ (see Figure 2).

\begin{center}
\setlength{\unitlength}{1mm}
\begin{picture}(30,35)(-15,-35)
\linethickness{1pt} \thinlines

\put(-40,0){\line(1,0){70}}
\put(-40,-40){\line(1,0){70}}
\put(-40,0){\line(0,-1){40}}
\put(30,0){\line(0,-1){40}}

\put(28,-2){\makebox(0,0)[c]{$G$}}
\put(-18,-17){\makebox(0,0)[c]{$a$}}
\put(-30,-30){\makebox(0,0)[c]{$A'$}}

\put(-15,-15){\circle*{2}}
\put(-15,-15){\circle{12}}

\put(-15,-15){\line(1,1){10}}
\put(-15,-15){\line(2,1){15}}
\put(-15,-15){\line(2,0){15}}
\put(-15,-15){\line(2,-1){20}}
\put(-15,-15){\line(2,-3){10}}
\put(-15,-15){\line(-1,-2){8}}
\put(-15,-15){\line(-1,0){18}}
\put(-15,-15){\line(0,-1){10}}
\put(-15,-15){\line(-1,1){12}}
\put(-15,-15){\line(1,5){2}}
\end{picture}
\end{center}

\vskip 0.1in

\begin{center}
Figure 2
\end{center}

\noindent It means that if $A'=A\backslash\{a\}$ and $G'=\langle A',E'\rangle$ is as before, then $G'$ can be again of any possible form, so as before
$$\frac{|{\rm Mod}_n(\phi_{\rm all}^{(1)})|}{|\textbf{S}_n(L)|}\leq\frac{n}{2^{n-1}}\longrightarrow_n 0.$$

\vskip 0.1in

(iii)  $\mu(\phi_{\rm one}^{(1)})=0$:  We have $\phi_{\rm one}^{(1)}:=(\exists x)(\exists! y)(xEy)$ and we have to  show that
\begin{equation} \label{E:showlimit3}
\frac{|{\rm Mod}_n(\phi_{\rm one}^{(1)})|}{|\textbf{S}_n(L)|}\longrightarrow_n 0.
\end{equation}
Let $G=\langle A,E\rangle$ be a graph, with $|A|=n$,  satisfying $\phi_{\rm one}^{(1)}$. Pick an  $a\in A$ witnessing $\phi_{\rm one}^{(1)}$. If  again $A'=A\backslash\{a\}$,  $a$ is adjacent to a unique element $b\in A'$ (see Figure 3), so there are $n-1$ choices for the given $a$, to be connected.

\begin{center}
\setlength{\unitlength}{1mm}
\begin{picture}(30,35)(-15,-35)
\linethickness{1pt} \thinlines

\put(-40,0){\line(1,0){70}}
\put(-40,-40){\line(1,0){70}}
\put(-40,0){\line(0,-1){40}}
\put(30,0){\line(0,-1){40}}

\put(28,-2){\makebox(0,0)[c]{$G$}}
\put(-18,-17){\makebox(0,0)[c]{$a$}}
\put(-25,-30){\makebox(0,0)[c]{$A'$}}
\put(8,-22){\makebox(0,0)[c]{$b$}}

\put(-15,-15){\circle*{2}}
\put(-15,-15){\circle{12}}
\put(-15,-15){\line(5,-1){20}}
\put(6,-19){\circle*{2}}

\end{picture}
\end{center}

\vskip 0.1in

\begin{center}
Figure 3
\end{center}

\noindent For each such choice the subgraph $G'=\langle A',E'\rangle$ can have any possible form, so the number of graphs on $A$ in which the specific $a$ witnesses $\phi_{\rm one}^{(1)}$ is at most $(n-1)\cdot2^{\binom{n-1}{2}}$. Letting again $a$ range over $A$, it follows that $|{\rm Mod}_n(\phi_{\rm one}^{(1)})|\leq n(n-1)2^{\binom{n-1}{2}}$. So
$$\frac{|{\rm Mod}_n(\phi_{\rm one}^{(1)})|}{|\textbf{S}_n(L)|}\leq\frac{n(n-1)2^{\binom{n-1}{2}}}
{2^{\binom{n}{2}}}=\frac{n(n-1)}{2^{n-1}}.$$
Since clearly $\lim_{n\rightarrow\infty}\frac{n(n-1)}{2^{n-1}}=0$, (\ref{E:showlimit3}) is true. \telos

\vskip 0.2in

Note that properties $\phi_{\rm none}(x)$ and $\phi_{\rm all}(x)$ can be shown to have typical degree $0$ also by a direct computation, without much effort, without the help of Lemma \ref{L:key1}. However such a  direct proof for $\phi_{\rm one}(x)$ seems  infeasible, because of the vast complexity of the required computations.  Moreover the method of proof for $\phi_{\rm one}(x)$ used in the last Proposition  can be easily generalized in order to apply to the  property ``$x$ is adjacent to exactly $k$ nodes'', for every $k\geq 1$.

\begin{Prop} \label{P:general}
For each $k\geq 1$, let $\phi_k(x)$ be the property  ``$x$ is adjacent to exactly $k$ nodes''. Then $d(\phi_k:\mbox{typ})=0$.
\end{Prop}

{\em Proof.} By Lemma \ref{L:key1} (i), it suffices to show that $\mu(\phi_k^{(1)})=0$. We argue as in the proof of \ref{P:application} concerning  $\phi_{\rm one}(x)$. Given a graph $G=\langle A,E\rangle$ with $|A|=n$ which satisfies $\phi_k^{(1)}$, let $a\in A$ witness $\phi_k^{(1)}$. Then $a$ is adjacent exactly to  $k$ elements of $A\backslash\{a\}$. Since there are $\binom{n-1}{k}$ $k$-element subsets of $A\backslash\{a\}$, there exist at most $\binom{n-1}{k}2^{\binom{n-1}{2}}$ graphs on $A$ in which $a$ witnesses $\phi_k^{(1)}$. Therefore  $|{\rm Mod}_n(\phi_k^{(1)})|\leq n\binom{n-1}{k}2^{\binom{n-1}{2}}$, and hence
$$\mu_n(\phi_k^{(1)})=\frac{|{\rm Mod}_n(\phi_k^{(1)})|}{|\textbf{S}_n(L)|}\leq \frac{n\binom{n-1}{k}2^{\binom{n-1}{2}}}
{2^{\binom{n}{2}}}=\frac{n\binom{n-1}{k}}{2^{n-1}}.$$
In $\frac{n\binom{n-1}{k}}{2^{n-1}}$ the nominator is a polynomial in $n$ of degree $k+1$, so by l' Hospital rule $\frac{n\binom{n-1}{k}}{2^{n-1}}\longrightarrow_n0$. \telos

\vskip 0.2in

\begin{Cor} \label{C:regular}
All properties $\phi_{\rm none}(x)$, $\phi_{\rm all}(x)$  and $\phi_k(x)$, for $k\geq 1$,  are regular.
\end{Cor}

{\em Proof.} We showed in the proofs of  Propositions \ref{P:application} and \ref{P:general} that for each of the aforementioned  properties $\phi(x)$,  $\mu(\phi^{(1)})=0$. This condition implies, according to \ref{L:key1} (ii), that  $\phi(x)$ is regular. \telos


\begin{thebibliography}{99}
\bibitem{Em95}
   H.-D. Ebbinghaus and J. Flum, {\em Finite Model Theory}, Springer 1995.
\bibitem{Fa76}
   R. Fagin, Probabilities on finite models, {\em J. Symb. Logic} \textbf{41} (1976), no. 1, 50 -58.
\bibitem{GETAL69}
   Y. V. Glebskii, D. I. Kogan, M. I. Liogon'kii, and V. A. Talanov, Range and degree of realizability of formulas in the restricted predicate calculus, {\em Cybernetics} \textbf{5} (1969), 142 -154.
\bibitem{Hetal08}
   J. M. Harris, J. L. List and M.J. Mossinghoff, {\em Combinatorics and Graph Theory}, Second Edition, Springer 2008.
\bibitem{Ru95}
   B. Russell, {\em My Philosophical Development,} Routledge, revised edition 1995.
\bibitem{Tz20}
   A. Tzouvaras, Russell's typicality as another randomness notion, {\em Math. Log. Quarterly} \textbf{66} (2020), no. 3, 355 -365.
\bibitem{Tz22}
   A. Tzouvaras, Typicality \`{a} la Russell in set theory, {\em Notre Dame J. Form. Logic} \textbf{63} (2022), no. 2, 185 -196.
\bibitem{Wiki}
   \verb"https://en.wikipedia.org/wiki/Central_binomial_coefficient".
\bibitem{Stirling}
    \verb"https://en.wikipedia.org/wiki/Stirling_numbers_of_the_second_kind".
\end{thebibliography}
\end{document}